\input amstex
\loadeufm

\documentstyle{amsppt}

\magnification=\magstep1
\baselineskip=18pt

\baselineskip=20pt
\parskip=5.5pt
\hsize=6.5truein
\vsize=9truein
\NoBlackBoxes

\topmatter

\title 
A Relationship Between the Dirichlet and
Regularity Problems for Elliptic Equations
\endtitle

\author Zhongwei Shen
\endauthor

\thanks
The author is supported in part by the NSF (DMS-0500257).
\endthanks

\address
Department of Mathematics, University of Kentucky,
Lexington, KY 40506
\endaddress

\email
shenz\@ms.uky.edu
\endemail

\leftheadtext{Zhongwei Shen} 
\rightheadtext{A Relationship Between the Dirichlet and Regularity Problems}

\abstract
Let $\Cal{L}=\text{div}A\nabla$ be a real, symmetric second order 
elliptic operator with bounded measurable coefficients.
Consider the elliptic equation
$\Cal{L}u=0$ in a bounded Lipschitz domain $\Omega$ of $\Bbb{R}^n$.
We study the relationship between the solvability of
the $L^p$ Dirichlet problem $(D)_p$
with boundary data in $L^p(\partial \Omega)$ and that of the $L^q$ regularity
problem $(R)_q$ with boundary data
in $W^{1,q}(\partial \Omega)$, where $1<p,q<\infty$. 
It is known that the solvability of $(R)_p$
implies that of $(D)_{p^\prime}$.
In this note we show that if
$(D)_{p^\prime}$ is solvable, then either $(R)_p$ is solvable or
$(R)_q$ is not solvable for any $1<q<\infty$.
\endabstract


\endtopmatter

\document

\centerline{\bf 1. Introduction}

Let $\Cal{L}=\text{div}A\nabla$ be a real, symmetric second
order elliptic operator with bounded measurable coefficients.
For $f\in C(\partial\Omega)$, consider the classical Dirichlet problem
$$
\left\{\aligned
& \Cal{L} u=0\ \ \ \text{ in }\Omega,\\
& u=f \ \ \ \text{ on } \partial\Omega,
\endaligned
\right.
\tag 1.1
$$
where $\Omega$ is a bounded Lipschitz domain in $\Bbb{R}^n$, $n\ge 2$.
Let $\| \cdot\|_p$ denote the norm in $L^p(\partial\Omega)$
with respect to the surface measure $d\sigma$ on $\partial\Omega$.
The $L^p$ Dirichlet problem $(D)_p$ is said to be solvable
if the unique solution $u\in C(\overline{\Omega})$ to (1.1)
satisfies the estimate $\| N(u)\|_p\le C\, \| f\|_p$.
Here $N(u)$ denotes the nontangential maximal function
of $u$.
If boundary data $f$ is in $W^{1,p}(\partial\Omega)$, i.e.,
the first order (tangential) derivatives of $f$ are also in $L^p(\partial
\Omega)$, it is natural to require that the solution satisfies the condition
$\|N(\nabla u)\|_p\le C\, \| f\|_{W^{1,p}(\partial\Omega)}$.
This is the so-called $L^p$ regularity problem.
We remark that since $\nabla u$ may not be locally bounded,
$N(\nabla u)$ needs to be suitably defined for weak solutions.
Both the $L^p$ Dirichlet and regularity problems (as well as the 
$L^p$ Neumann
problem) have been studied extensively. We refer the reader to
the monograph \cite{\bf K} by Kenig for a survey of results 
as well as a list of open problems.

The purpose of this note is to study the relationship between the Dirichlet
problem and the regularity problem. 
It is well known that the solvability of the $L^p$ regularity problem
implies that of the $L^{p^\prime}$ Dirichlet problem (see e.g. \cite{\bf KP2}),
where $p^\prime=p/(p-1)$ denotes the dual exponent of $p$.
The converse is also true in the case of
 Laplace's equation $\Delta u=0$ on Lipschitz domains
\cite{\bf V}.
However whether the converse holds for a general second order
elliptic operator with nonsmooth coefficients remains open.
In this note we obtain a partial result.

\proclaim{\bf Main Theorem} 
Let $\Cal{L}$ be a real, symmetric second order elliptic operator
of divergence form with bounded measurable coefficients.
Let $1<p<\infty$ and $\Omega$ be a bounded Lipschitz domain.
Suppose that the $L^{p^\prime}$ Dirichlet problem for $\Cal{L}$
on $\Omega$ is solvable.
Then either the $L^p$ regularity problem $(R)_p$
 is solvable or $(R)_q$ is not solvable for any $1<q<\infty$.
\endproclaim

We remark that for Laplace's equation in a Lipschitz domain $\Omega$, the
Dirichlet problem $(D)_{p^\prime}$ and the regularity problem
$(R)_p$ are solvable for $1<p<2+\varepsilon$, where
$\varepsilon>0$ depends on $\Omega$ \cite{\bf D1, D2, JK, V}.
If $\Omega$ is a $C^1$ domain, $(D)_p$ and $(R)_p$ 
are solvable for any $1<p<\infty$ \cite{\bf FJR}.
However, for a general second order elliptic operator
$\Cal{L}$,
$(D)_p$ (thus $(R)_{p^\prime}$)
may not be solvable for any $1<p<\infty$,
even if the coefficients of $\Cal{L}$ are continuous and
$\Omega$ is smooth.
Furthermore, it is known that the Dirichlet problem
$(D)_{p^\prime}$ for $\Cal{L}$ on $\Omega$
is solvable if and only if
the $\Cal{L}$-harmonic measure 
is a $B_p$ weight with respect to the surface measure
on $\partial\Omega$. We
refer the reader to \cite{\bf K}
for references on these and other deep results on the solvability of 
the $L^p$ Dirichlet problem. Concerning
the $L^p$ regularity problem for
general second order elliptic operators, we mention that the study
was initiated by Kenig and Pipher in \cite{\bf KP1, KP2}. 
In particular it was proved in \cite{\bf KP1} that
the solvability of $(R)_p$ implies that
of $(R)_q$ for all $1<q<p+\varepsilon$.
This fact is used in the proof of the main theorem.
 
Our main theorem will be proved in
two steps.
First we  establish a weak reverse
H\"older estimate,
$$
\left(\frac{1}{|I_r|}
\int_{I_r}
|N(\nabla u)|^p\, d\sigma\right)^{1/p}
\le \frac{C}{|I_{6r}|}
\int_{I_{6r}} |N(\nabla u)|\, d\sigma,
\tag 1.2
$$ 
where $I_r$ is a surface cube on $\partial\Omega$, and
$u$ is a weak solution of $\Cal{L}u=0$ in $\Omega$ whose
boundary data vanishes on $I_{6r}$.
This is done by using the comparison principle
and properties of $\Cal{L}$-harmonic measures (see Theorem 2.9).
It worths pointing out that for (1.2) to hold, one only needs to assume
that $(D)_{p^\prime}$ is solvable.

The second step of the proof of the main theorem relies on a real
variable argument.
It uses a maximal dyadic cube decomposition on $\partial\Omega$
and the reverse
H\"older estimate (1.2) to establish a good-$\lambda$
type inequality (see Lemma 3.4).
It is here that one needs to assume the 
solvability of $(R)_q$ for some $q<p$. 
Motivated by \cite{\bf CP}
(also see \cite{\bf W}),
this approach of combining the Calder\'on-Zygmund decomposition
with the reserve H\"older estimates was developed in \cite{\bf S1, S2, S3}
to study the $L^p$ bounds of Riesz transforms associated with
operator $\Cal{L}$, and the $L^p$ estimates
for elliptic systems and higher order elliptic equations
on Lipschitz domains.
We remark that a similar method was used simultaneously 
and independently with different motivation by
Auscher, Coulhon, Duong and Hofmann \cite{\bf ACDH, A}
in the study of Riesz transforms on manifold as well as
elliptic operators with complex coefficients.

\medskip
 
\centerline{\bf 2. A weak reverse H\"older inequality}

Throughout this note, we will use
$\Omega$ to denote a bounded Lipschitz domain in $\Bbb{R}^n$.
 For $P\in \partial
\Omega$, let
$$
\Gamma_\alpha (P)=\left\{ x\in \Omega:\ \ |x-P|<(1+\alpha)\delta(x)\, 
\right\},
\tag 2.1
$$
where $\alpha>0$ and $\delta(x)=\text{dist}(x,\partial \Omega)$.
We define 
$$
N_\alpha (u)(P)=\sup\left\{
\left(\frac{1}{|B(x,\frac12 \delta(x))|}
\int_{B(x,\frac12 \delta(x))}
|u(y)|^2\, dy\right)^{1/2}:\ \ 
x\in \Gamma_\alpha(P)\, \right\},
\tag 2.2
$$
for any $u\in L^2_{\text{loc}}(\Omega)$.
We will use $N(u)$ for $N_1(u)$.
This is a variant of the usual nontangential
maximal function $(u)^*$, which is defined
by
$$
(u)^*(P)=\sup\big\{ |u(x)|:\ x\in \Gamma_1 (P)
\, \big\}.
$$
It is known that if $\Cal{L} u=0$ in $\Omega$,
then $\|N_\alpha (u)\|_p$ and $\|(u)^*\|_p$ are equivalent 
for any $ 0<p\le \infty$.

Assume $0\in \partial\Omega$ and
$$
\Omega\cap B(0,r_0)=\big\{ (x^\prime,x_n)\in \Bbb{R}^n:\
x_n>\psi(x^\prime)\, \big\}
\cap B(0,r_0),
\tag 2.3
$$
where $B(0,r_0)$ denotes the ball centered at $0$ with radius $r_0$
and $\psi:\Bbb{R}^{n-1}\to \Bbb{R}$ is a Lipschitz function, $\psi(0)=0$.
For $r>0$, we let
$$
\aligned
&I_r=\big\{ (x^\prime, \psi(x^\prime)):\
|x_1|<r, \dots, |x_{n-1}|<r\,\big\},\\
&D_r=\big\{ (x^\prime, x_n)\in \Bbb{R}^n:\
|x_1|<r, \dots, |x_{n-1}|<r,\
\psi(x^\prime)<x_n<\psi(x^\prime)+r\, \big\}.
\endaligned
\tag 2.4
$$
Note that if $0<r<c\, r_0$, then 
$D_r\subset \Omega$ and $\partial D_r\cap \partial\Omega
=I_r$.

\proclaim{\bf Lemma 2.5}
Let $u\in C(\overline{\Omega})$ be a weak solution of $\Cal{L}
u=0$ in $\Omega$. Suppose that $u=0$ on $I_{5r}$ for some
$0<r<c\, r_0$. Then for any $x\in D_{2r}$,
$$
|u(x)|\le C \cdot
\frac{G(x,z)}{G(A_r, z)}\left\{ \frac{1}{r^n}\int_{D_{4r}}
|u(y)|^2\, dy\right\}^{1/2},
\tag 2.6
$$
where $A_r=(0,\frac{r}{2})$, $z\in \Omega\setminus
D_{10r}$, and $G(x,y)$ denotes the Green function
for $\Cal{L}$ on $\Omega$.
\endproclaim

\demo{Proof}
Write $u=u_1-u_2$ on $D_{3r}$,
where $u_1,u_2$ are nonnegative weak solutions
on Lipschitz domain $D_{3r}$
 with boundary values $u_+=\max( u, 0)$ and $u_-
=\max (-u, 0)$ on 
$\partial D_{3r}$ respectively.
By the comparison principle for nonnegative weak
solutions \cite{\bf CFMS}, we have
$$
u_j(x)\le C\cdot \frac{G(x,z)}{G(A_r,z)}\cdot u_j(A_r)
\le C \cdot\frac{G(x,z)}{G(A_r,z)}\cdot\max_{\partial D_{3r}}
|u_j|, \ \ \ \ j=1,2,
\tag 2.7
$$
where $x\in D_{2r}$ and  $z\in \Omega\setminus D_{10r}$.
It follows that
$$
\aligned
|u(x)| &\le 
|u_1(x)| +|u_2(x)|
\le C\cdot
\frac{G(x,z)}{G(A_r,z)}\cdot \max_{\partial D_{3r}} |u|\\
&\le C\cdot \frac{G(x,z)}{G(A_r,z)}
\left\{ \frac{1}{r^n}\int_{D_{4r}}
|u(y)|^2\, dy\right\}^{1/2},
\endaligned
\tag 2.8
$$ 
where we have used the boundary $L^\infty$ estimate 
in the last step.
The proof is finished.
\enddemo

\proclaim{\bf Theorem 2.9} 
Suppose $1<p<\infty$ and 
$(D)_{p^\prime}$ is solvable for operator $\Cal{L}$ on $\Omega$.
Let $u\in C(\overline{\Omega})$ be a weak
solution of $\Cal{L}u=0$ in $\Omega$.
If $u=0$ on $I_{6r}$ for some $0<r<c\, r_0$, then
$$
\left(\frac{1}{|I_r|}\int_{I_r}
|N(\nabla u)|^p\, d\sigma \right)^{1/p}
\le \frac{C}{|I_{6r}|}
\int_{I_{6r}}
|N(\nabla u)|\, d\sigma,
\tag 2.10
$$
where $|I_r|$ denotes the surface measure of $I_r$.
\endproclaim
 
\demo{Proof}
We begin by observing that for any $P\in I_r$,
$$
N(\nabla u)(P)
\le C\, \left\{ \left(\frac{u}{\delta}\right)^*_{20,r}(P)
+\frac{1}{|I_{6r}|}
\int_{I_{6r}} N(\nabla u)\, d\sigma\right\},
\tag 2.11
$$
where $(\frac{u}{\delta})^*_{20,r}(P)
=\sup\left\{ \frac{|u(x)|}{\delta(x)}:\
|x-P|\le c\, r \text{ and }
x\in \Gamma_{20}(P)\, \right\}$.
To see (2.11), we note that if $x\in \Gamma_1(P)$ and $|x-P|\ge c\, r$
for some $P\in I_r$, then
$$
\left(\frac{1}{|B(x,\frac12 \delta(x))|}
\int_{B(x,\frac12 \delta(x))}
|\nabla u(y)|^2\, dy\right)^{1/2}
\le \frac{C}{|I_{6r}|}
\int_{I_{6r}} N(\nabla u)\, d\sigma.
\tag 2.12
$$
Indeed, if the left side of (2.12) is greater than $\lambda$, then
$|\{ P\in I_{6r}:\ N(\nabla u)(P)>c\, \lambda\}\ge c\, |I_{6r}|$.
This may be seen by subdividing $B(x,\frac12 \delta(x))$.
On the other hand, if $x\in \Gamma_1(P)$ and $|x-P|< c\, r$ for some
$P\in I_r$, the left side of (2.12) is bounded by
$$
 C\, \left(\frac{1}{|B(x,\frac34 \delta(x))|}
\int_{B(x,\frac34 \delta(x))}
\left|\frac{u(y)}{\delta(y)}\right|^2\, dy\right)^{1/2}
\le C \, \left(\frac{u}{\delta}\right)^*_{20,r}
(P).
\tag 2.13
$$
This follows from the Cacciopoli
inequality.

To estimate $(\frac{u}{\delta})^*_{20,r}$ on $I_r$, we use
Lemma 2.5.
This gives us
$$
\left(\frac{u}{\delta}\right)^*_{20,r}
(P)\le \frac{C}{G(A_r,z)}
\left(\frac{1}{r^n}\int_{D_{4r}} |u(y)|^2\, dy\right)^{1/2}
\cdot \left(\frac{G(\cdot,z)}{\delta(\cdot)}\right)^*_{20,r}
(P)
\tag 2.14
$$
for any $P\in I_r$. Now, fix $z
\in \Omega\setminus D_{10r}$.
Let $\omega=\omega^z$ denote the $\Cal{L}$-harmonic
measure on $\partial\Omega$, evaluated at $z$. Also, for $x=(x^\prime, x_n)
\in D_{5r}$, let $\hat{x}=(x^\prime, \psi(x^\prime))
\in \partial\Omega$
and $S(\hat{x},t)=B(\hat{x},t)\cap \partial\Omega$.
Since
$$
\frac{G(z,x)}{\delta(x)}
\approx \frac{\omega(S(\hat{x},\delta(x))}{\delta(x)^{n-1}},
\tag 2.15
$$
for $x\in D_{2r}$ and $\omega$ is a doubling measure \cite{\bf CMFS}, we have
$$
\left(\frac{G(z,\cdot)}{\delta(\cdot)}\right)^*_{20,r}(P)
\le C\, M_{\sigma,r} (\omega)(P),
\tag 2.16
$$
where 
$$
M_{\sigma,r}(\omega)(P)=\sup\left\{ 
\frac{\omega^z(S(P,t))}{t^{n-1}}:\ \ 0<t<r\, \right\}
\tag 2.17
$$
is a localized Hardy-Littlewood maximal function of $\omega$ with respect
to the surface measure.

Finally, since $(D)_{p^\prime}$ is solvable for $\Cal{L}$
on $\Omega$ , $\omega$ satisfies
the Reverse H\"older inequality,
$$
\left(\frac{1}{|S|}\int_S \big|\frac{d\omega}{d\sigma}\big|^p\,
d\sigma\right)^{1/p}
\le \frac{C\, \omega(S)}{|S|},
\tag 2.18
$$
for any surface ball $S=B(P,t)\cap\partial\Omega$ on $\partial\Omega$.
This, together with (2.14)-(2.16) as well as the $L^p$ boundedness
of the Hardy-Littlewood maximal operator, yields
$$
\left(\frac{1}{|I_{r}|}
\int_{I_{r}}
\big|\left(\frac{u}{\delta}\right)^*_{20,r}\big|^p\, d\sigma\right)^{1/p}
\le \frac{C}{r}
\left(\frac{1}{r^n}\int_{D_{4r}}
|u(y)|^2\, dy\right)^{1/2}.
\tag 2.19
$$
To finish the proof, we note that
$$
\aligned
& \frac{1}{r}\left(\frac{1}{r^{n}}\int_{D_{4r}}
|u|^2\, dy\right)^{1/2}
\le C\, \left(\frac{1}{r^n}\int_{D_{4r}}
|\nabla u|^2\, dy\right)^{1/2}\\
&\le \frac{C}{r^n}\int_{D_{5r}} |\nabla u|\, dy
\le \frac{C}{r^{n-1}}\int_{I_{6r}}
N(\nabla u)\, d\sigma.
\endaligned
\tag 2.20
$$
We point out that
 we have used Poincar\'e inequality for the first inequality in (2.20)
and some well 
known properties of weak solutions for the second. Also, 
 the third inequality in (2.20) follows from
$$
\int_{I_{6r}} N(\nabla u)\, d\sigma
\ge \frac{C}{r}
\int_{D_{6r}} \left\{ \frac{1}{\delta(x)^n}\int_{B(x,\frac12 \delta(x))}
|\nabla u(y)|\, dy\right\}\, dx
\tag 2.21
$$
and the Fubini's theorem.
The desired reverse H\"older estimate (2.10) now follows
by combining
(2.11), (2.19) and (2.20).
\enddemo
\medskip

\centerline{\bf 3. The proof of the Main Theorem}

Let $1<p<\infty$. Suppose that the Dirichlet problem
 $(D)_{p^\prime}$ is solvable for 
operator $\Cal{L}$ on Lipschitz domain $\Omega$. Also assume that
the regularity problem $(R)_q$ is solvable for some
$1<q<p$. We will show that $(R)_p$ is solvable.
Let $f\in W^{1,p}(\partial\Omega)\cap
C(\partial\Omega)$ and $u\in C(\overline{\Omega})$
be the unique weak solution of 
the classical Dirichlet problem (1.1). We need to prove
that 
$$
\|N(\nabla u)\|_p\le C\, \|\nabla_t f\|_p,
\tag 3.1
$$
where $\nabla_t f$ denotes the first order tangential derivatives of 
$f$.

To this end, we
fix $P_0\in \partial\Omega$. By translation and rotation, we may
assume that $P_0$ is the origin and (2.3) holds.
Let $D=\{ (x^\prime,x_n)\in \Bbb{R}^n:\ x_n>\psi(x^\prime)\}$
and $D_r$, $I_r$ be defined as in (2.4).
Define the map $\Phi:\partial D\to \Bbb{R}^{n-1}$ by
$\Phi(x^\prime,\psi(x^\prime))=x^\prime$.
We say that $Q\subset\partial D$ is a ``cube'' of $\partial D$ if
$\Phi(Q)$ is a cube of $\Bbb{R}^{n-1}$.
The dyadic subcubes of a cube on $\partial D$ are defined similarly.
We will use $\rho Q$ for $\Phi^{-1}[ \rho \Phi(Q)]$.
If $Q$ is a cube on $\partial D$, we define a localized
Hardy-Littlewood maximal function $M_Q$ by
$$
M_Q(f)(P)=\sup\Sb P\in Q^\prime\\ \text{cube } Q^\prime\subset Q\endSb
\frac{1}{|Q^\prime|} \int_{Q^\prime}
|f|\, d\sigma.
\tag 3.2
$$
For $\lambda>0$ and $0<r<c\, r_0$, let
$$
E(\lambda)=\left\{ P\in I_r:\
M_{{I_{2r}}}
(|N(\nabla u)|^q)(P)>\lambda\, \right\}.
\tag 3.3
$$

\proclaim{\bf Lemma 3.4}
Let $1<q<p<\infty$. Suppose that $(D)_{p^\prime}$ and $(R)_q$
are solvable.
There exist positive constants $\varepsilon$, $\gamma$, $C_0$
depending only on $p$, $q$, $n$, $\Cal{L}$ and $\Omega$ such that
$$
|E(A\lambda)|
\le \varepsilon |E(\lambda)|
+|\left\{ P\in I_r :\
M_{I_{2r}}(|\nabla_t f|^q)(P)>\gamma \lambda\, 
\right\}|
\tag 3.5
$$
for all $\lambda\ge \lambda_0$, where $A=1/(2\varepsilon)^{q/p}$ and
$$
\lambda_0
=\frac{C_0}{|I_{2r}|}
\int_{I_{2r}} |N(\nabla u)|^q\, d\sigma.
\tag 3.6
$$
\endproclaim

\demo{Proof}
Let $\varepsilon\in (0,1)$ be a small constant to be determined.
By the weak $(1,1)$ estimate of the Hardy-Littlewood
maximal function, we have
$$
|E(\lambda)|\le \frac{C}{\lambda}
\int_{I_{2r}} |N(\nabla u)|^q\, d\sigma,
\tag 3.7
$$
where $C$ depends only on $n$ and $\|\nabla \psi\|_\infty$.
It follows that $|E(\lambda)|<\varepsilon |I_r|$
if $\lambda\ge \lambda_0$, where $\lambda_0$ is given by (3.6)
with a large $C_0$.

Next we perform a Calder\'on-Zygmund decomposition
on $E(\lambda)$ as a relative open subset of $I_r$.
 This produces a collection of
disjoint dyadic subcubes $\{ Q_k\}$ of $I_r$ such that
$E(\lambda)=\bigcup_k Q_k$ and each $Q_k$ is maximal.
We may assume that $\varepsilon$ is sufficiently small
so that $32 Q_k\subset I_{2r}$.
We claim that it is possible to choose constants $\varepsilon$, 
$\gamma$, $C_0$ so that
if $Q_k$ is a cube with the property 
$$\{ P\in Q_k: \ M_{I_{2r}}(|\nabla_t f|^q)(P)\le \gamma\lambda\, \}
\neq \emptyset,
\tag 3.8
$$
 then
$$
|E(A\lambda)\cap Q_k|\le \varepsilon |Q_k|.
\tag 3.9
$$
From this, estimate (3.5) follows by summation.

To establish (3.9), we first observe that
$$
M_{I_{2r}}(|N(\nabla u)|^q)(P)\le \max
\left(M_{2Q_k}(|N(\nabla u)|^q)(P), C_1\lambda\right)
\tag 3.10
$$
for any $P\in Q_k$, where $C_1$ depends only on
$n$ and $\|\nabla \psi\|_\infty$.
 This is because $Q_k$ is maximal and so $3Q_k\nsubseteq E(\lambda)$.
Assume $A\ge C_1$.
It follows that
$$
|Q_k\cap E(A\lambda)|
\le |\left\{ P\in Q_k:\ M_{2Q_k} (|N(\nabla u)|^q)(P)
>A\lambda\, \right\}|.
\tag 3.11
$$
Now let $v=v_k$ be the unique weak solution of $\Cal{L}v=0$
in $\Omega$ with boundary data $\varphi (f-\alpha)$, where
$$
\alpha =\frac{1}{|17Q_k|}
\int_{17 Q_k} f\, d\sigma,
\tag 3.12
$$
and $\varphi=\varphi_k$ is a smooth cut-off function on
$\Bbb{R}^n$ such that
$\varphi=1$ on $16 Q_k$, $\varphi=0$ on $\partial\Omega
\setminus 17 Q_k$, and $|\nabla \varphi |\le C/\ell_k$, where
$\ell_k=|Q_k|^{1/(n-1)}$.
Let $\bar{p}>p$. In view of (3.11), we have
$$
\aligned
&|Q_k\cap E(A\lambda)|
 \le |\left\{ P\in Q_k:\
M_{2Q_k}(|N(\nabla u-\nabla v)|^q)(P)>\frac{A\lambda}
{2^q}\, \right\}|\\
&\ \ \ \ 
+|\left\{ P\in Q_k:\
M_{2Q_k}(|N(\nabla v)|^q)(P)
>\frac{A\lambda}{2^q}\, \right\}|\\
&\le \frac{C}{(A\lambda)^{\bar{p}/q}}
\int_{2Q_k}
|N(\nabla u-\nabla v)|^{\bar{p}}\, d\sigma
+\frac{C}{A\lambda}
\int_{2Q_k}|N(\nabla v)|^q\, d\sigma,
\endaligned
\tag 3.13
$$
where we have used the weak $(\frac{\bar{p}}{q},\frac{\bar{p}}{q})$
and $(1,1)$ estimates for the Hardy-Littlewood
maximal operator. Since $(R)_q$ is solvable, the second term in the right side
of (3.13) is bounded by
$$
\aligned
\frac{C}{A\lambda}\int_{\partial\Omega}
|\nabla_t [\varphi(f-\alpha)]|^q\, d\sigma
&\le \frac{C}{A\lambda}\int_{17Q_k} |\nabla_t f|^q\, d\sigma\\
&\le \frac{C}{A\lambda}\cdot \gamma \lambda |17Q_k|
\le \frac{C\gamma}{A}\cdot |Q_k|,
\endaligned
\tag 3.14
$$
where we have used Poincar\'e inequality for the first
inequality and (3.8) for the second.

To handle the first term on the right side of (3.13), 
we observe that $u-v-\alpha$ is a weak solution whose boundary data
$(f-\alpha)(1-\varphi)$ vanishes
on $16 Q_k$. Also note that the solvability of $(D)_{p^\prime}$
implies that $(D)_{\bar{p}^\prime}$ for some $\bar{p}>p$.
It then follows by Theorem 2.9 that  the first term on the right
side of (3.13) is bounded by
$$
\aligned
&\frac{C|Q_k|}{(A\lambda)^{\bar{p}/q}}
\left(\frac{1}{|12Q_k|}
\int_{12Q_k} |N(\nabla u-\nabla v)|\, d\sigma\right)^{\bar{p}}\\
&\le 
\frac{C|Q_k|}{(A\lambda)^{\bar{p}/q}}
\left\{
\left(\frac{1}{|12 Q_k|}
\int_{12Q_k}|N(\nabla u)|^q\, d\sigma\right)^{\bar{p}/q}
+\left(\frac{1}{|12Q_k|}
\int_{12 Q_k} |N(\nabla v)|^q\, d\sigma\right)^{\bar{p}/{q}}
\right\}\\
&\le \frac{C|Q_k|}{(A\lambda)^{\bar{p}/q}}\cdot
\left\{ \lambda^{\bar{p}/{q}} +(\gamma \lambda)^{\bar{p}/q}\right\}
\le \frac{C}{A^{\bar{p}/q}}\cdot |Q_k|, 
\endaligned
$$
where in the second inequality,
we have used the solvability of $(R)_q$,  (3.14) as well as 
the fact that $Q_k$ is maximal.
This, together with (3.13) and (3.14), gives
$$
\aligned
|Q_k\cap E(A\lambda)|
&
\le |Q_k|\left\{ \frac{C\gamma}{A} +\frac{C}{A^{\bar{p}/q}}\right\}\\
&= \varepsilon|Q_k|\left\{
C_2\gamma \varepsilon^{-\frac{q}{p}-1}
+C_2\, \varepsilon^{\frac{\bar{p}}{p}-1}\right\},
\endaligned
\tag 3.16
$$
since $A=1/(2\varepsilon)^{q/p}$.

Finally, since $\bar{p}>p$, we may choose $\varepsilon>0$ so small
that $C_2\, \varepsilon^{\frac{\bar{p}}{p}-1}<1/2$. With this
$\varepsilon$ fixed, we then choose $\gamma>0$ so small that
$C_2\gamma \varepsilon^{-\frac{q}{p}-1}<1/2$.
The desired estimate (3.9) follows.
This completes the proof.
\enddemo

\demo{\bf Proof of The Main Theorem}
We multiply both sides of (3.5) by $\lambda^{\frac{p}{q}-1}$ and integrate
the resulting inequality in $\lambda\in (\lambda_0, \Lambda)$. 
This gives
$$
\int_{\lambda_0}^\Lambda |E(A\lambda)|\lambda^{\frac{p}{q}-1}\, d\lambda
\le \varepsilon
\int_{\lambda_0}^\Lambda
|E(\lambda)|\lambda^{\frac{p}{q}-1}\, d\lambda
+ C\int_{I_{2r}}
|\nabla_t f|^p\, d\sigma.
\tag 3.17
$$
Since $A^{p/q}=1/(2\varepsilon)$, by a change of variable, we obtain
$$
\int_0^{\Lambda}
|E(\lambda)|\lambda^{\frac{p}{q}-1}\, d\lambda
\le C\, \lambda_0^{\frac{p}{q}}\, |I_{2r}|
+C\, \int_{I_{2r}} |\nabla_t f|^p\, d\sigma.
\tag 3.18
$$
It follows by letting $\Lambda\to\infty$ in (3.18) that
$$
\int_{I_r} |N(\nabla u)|^p\, d\sigma
\le C\, \lambda_0^{\frac{p}{q}} |I_{2r}|
+C\int_{I_{2r}} |\nabla_t f|^p\, d\sigma.
\tag 3.19
$$
In view of (3.6), we have proved that
$$
\aligned
&\left(\frac{1}{|I_r|}
\int_{I_r}
|N(\nabla u)|^p\, d\sigma\right)^{1/p}\\
&\le C\, \left\{\left(\frac{1}{|I_{2r}|}
\int_{I_{2r}}
|N(\nabla u)|^q\, d\sigma\right)^{1/q}
+\left(\frac{1}{|I_{2r}|}
\int_{I_{2r}}
|\nabla_t f|^p\, d\sigma\right)^{1/p}
\right\}.
\endaligned
\tag 3.20
$$
From this and the estimate
$\| N(\nabla u)\|_q
\le C\, \|\nabla_t f\|_q$,
 we obtain $\|N(\nabla u)\|_p\le C\, \|\nabla_t f\|_p$
by covering $\partial\Omega$ with a finite number
of coordinate patches.
The proof is finished.
\enddemo

\enddocument

\end